\documentclass[11pt]{article}

\usepackage[margin=1.15in]{geometry}
\usepackage{amsmath,amssymb,amsthm,mathtools,mathrsfs}
\usepackage{enumitem}
\usepackage{microtype}
\usepackage[hidelinks]{hyperref}
\usepackage[nameinlink,capitalise]{cleveref}

\newtheorem{theorem}{Theorem}[section]
\newtheorem{proposition}[theorem]{Proposition}
\newtheorem{lemma}[theorem]{Lemma}
\newtheorem{corollary}[theorem]{Corollary}
\newtheorem{remark}[theorem]{Remark}

\theoremstyle{definition}

\newtheorem{example}[theorem]{Example}

\newcommand{\Add}{\operatorname{Add}}
\newcommand{\gen}{\operatorname{gen}}
\newcommand{\Proj}{\mathsf{Proj}}
\newcommand{\PProj}{\mathsf{PProj}}
\newcommand{\SSmod}{\mathsf{SS}}
\newcommand{\Spec}{\operatorname{Spec}_{\mathrm{cat}}}

\newcommand{\abs}[1]{\lvert #1\rvert}

\title{Characterizing categoricity in the class  $Add(M)$}
\author{Xiaolei Zhang\\
\small{School of Mathematics and Statistics, Tianshui Normal University
Tianshui 741001, China}\\
\small{zxlrghj@163.com}}
\date{}

\begin{document}
\maketitle

\begin{abstract}
We study categoricity of the additive closure $\operatorname{Add}(M)$, consisting of all direct summands of arbitrary direct sums of copies of a fixed module $M$. For an $\eta$-generated module $M$, we prove that categoricity in a single cardinal $\lambda \ge \theta = \max\{|R|,\eta\}$ is equivalent to the stabilization condition $P^{(\eta)} \cong M^{(\eta)}$ for every nonzero $\eta$-generated $P \in \operatorname{Add}(M)$, or equivalently $\operatorname{Add}(P)=\operatorname{Add}(M)$. This condition implies that every $X \in \operatorname{Add}(M)$ of cardinality $\lambda > \nu_M = \max\{\theta, |\mathcal{S}_\eta(M)|\}$ is isomorphic to $M^{(\lambda)}$, where $\mathcal{S}_\eta(M)$ is a representative set of the nonzero $\eta$-generated modules in $\operatorname{Add}(M)$. Since $|\mathcal{S}_\eta(M)| \le 2^\theta$, we obtain a uniform ZFC tail above $2^\theta$. For projective modules, the criterion reduces to Bass's right $p$-connectedness. For pure-projective modules, it forces von Neumann regularity and simplicity: the pure-projective modules are eventually categorical if and only if $R$ is a simple von Neumann regular ring. In the commutative case, this recovers precisely the fields. Our approach combines Walker's decomposition theorem, support reduction, and Eilenberg absorption, and we highlight several subtle pitfalls in naive downward transfer arguments.
\end{abstract}


\section{Introduction}

Categoricity asks when all objects of a fixed cardinality in a given class are isomorphic.  Its model-theoretic origin goes back to \L o\'s and Morley, and its modern form is closely connected with Shelah's categoricity program \cite{Los,Morley,Shelah}.  In module theory, categoricity has recently been investigated for several homologically defined classes.  Mazari-Armida studied, among others, flat, absolutely pure, and locally pure-injective modules \cite{Mazari}.  Trlifaj gave a sharp algebraic criterion for categoricity of projective modules and proved transfer results for deconstructible classes \cite{Trlifaj}; the deconstructible setting was subsequently developed further by \v{S}aroch and Trlifaj \cite{SarochTrlifaj}.

The class studied in this paper has a different form.  For a right $R$-module $M$, let
\[
  \Add(M)=\{X\mid X\text{ is isomorphic to a direct summand of }M^{(I)}
  \text{ for some set }I\}.
\]
The class $\Add(M)$ is closed under arbitrary direct sums and direct summands, but it need not be closed under extensions and therefore need not be deconstructible.  Nevertheless, Walker's decomposition theorem \cite{Walker} gives a strong and very useful decomposition principle for direct summands of direct sums of boundedly generated modules.

Suppose that $M$ is generated by at most an infinite cardinal $\eta$ and put
\[
   \theta=\max\{|R|,\eta\}.
\]
Every object of $\Add(M)$ is then a direct sum of $\eta$-generated modules of cardinality at most $\theta$.  This immediately yields an \emph{upward} transfer: if $\Add(M)$ is categorical at one cardinal $\lambda_0\geq\theta$, then it is categorical at every cardinal $\lambda\geq\lambda_0$.  The delicate point is the opposite direction.  Categoricity at $\lambda_0$ cannot simply be applied to modules of cardinality $\lambda<\lambda_0$.  A genuine downward transfer must extract additional algebraic information from the isomorphisms occurring at $\lambda_0$.

The key invariant is the set of small direct-summand types.  Let $\mathscr S_\eta(M)$ be a representative set of the isomorphism classes of all nonzero $\eta$-generated modules in $\Add(M)$, and put
\[
  \sigma_M=|\mathscr S_\eta(M)|,
  \qquad
  \nu_M=\max\{\theta,\sigma_M\}.
\]
One-cardinal categoricity forces the stabilization identity
\[
  P^{(\eta)}\cong M^{(\eta)}
  \qquad(P\in\mathscr S_\eta(M)).
\]
Equivalently, every nonzero small object $P$ already generates the whole additive class:
\[
  \Add(P)=\Add(M).
\]
This local statement is the decisive algebraic content of categoricity.  Once it holds, every sufficiently large object can be written as a bounded ``remainder'' plus a copy of $M^{(\lambda)}$; Eilenberg absorption then removes the remainder.  We obtain categoricity in every cardinal $\lambda>\nu_M$.  Since $\sigma_M\leq2^\theta$, categoricity in one cardinal $\geq\theta$ implies categoricity in every cardinal $>2^\theta$.  If $\sigma_M\leq\theta$, then the natural threshold $\theta$ is recovered.

The projective case admits a particularly clean interpretation.  Bass called a ring $p$-connected when the trace ideal of every nonzero projective module is the whole ring \cite{Bass}; equivalently, every nonzero projective module is a generator.  We prove that the class $\Proj$ of projective right modules is categorical in one sufficiently large cardinal if and only if $R$ is right $p$-connected.  This recasts the stabilization criterion of Trlifaj \cite[Proposition~1.1]{Trlifaj} in intrinsic ring-theoretic terms.

The pure-projective case becomes even more rigid.  A right module is pure-projective precisely when it is a direct summand of a direct sum of finitely presented modules \cite{Warfield,GobelTrlifaj}.  We prove that eventual categoricity of pure-projective modules is equivalent to the conjunction of von Neumann regularity and right $p$-connectedness.  For von Neumann regular rings, right $p$-connectedness is in turn equivalent to simplicity.  Hence
\[
\begin{split}
 \PProj\text{ is categorical in one sufficiently large cardinal}
 \quad\Longleftrightarrow\quad{}&\\
 R\text{ is a simple von Neumann regular ring}.
\end{split}
\]
This gives a complete ring-theoretic answer, including non-artinian examples.  In the commutative case, the condition reduces to $R$ being a field.

We also treat semisimple modules directly.  Their class is categorical in one, equivalently every, cardinal strictly above $\max\{|R|,\aleph_0\}$ exactly when $R$ has a unique simple right module up to isomorphism.  The direct proof is important: the direct sum of representatives of all simple modules need not be generated by $\leq|R|+\aleph_0$ elements, so a bounded-generator theorem cannot be invoked without an additional cardinality argument.

The paper is organized as follows.  Section~2 records the decomposition and absorption tools.  Sections~3 and~4 prove the upward transfer and the small-summand stabilization theorem.  Section~5 identifies projective categoricity with right $p$-connectedness and treats classes $\Add(M)$ containing $R$.  Section~6 gives the semisimple characterization.  Section~7 proves the simple-von-Neumann-regular characterization for pure-projective modules and gives examples.  
\section{Preliminaries}

Throughout, $R$ is an associative ring with identity and all modules are unitary right $R$-modules.  For a module $X$, $\abs X$ denotes the cardinality of its underlying set, and $\gen(X)$ denotes its minimal number of generators.  If $I$ is a set, $X^{(I)}$ denotes the direct sum of $\abs I$ copies of $X$; for a cardinal $\kappa$ we also write $X^{(\kappa)}$.

A class $\mathcal C$ of modules is \emph{$\lambda$-categorical} if it contains a module of cardinality $\lambda$ and any two modules in $\mathcal C$ of cardinality $\lambda$ are isomorphic.  Its categoricity spectrum is
\[
  \Spec(\mathcal C)=\{\lambda\mid \mathcal C\text{ is }\lambda\text{-categorical}\}.
\]
A module is called \emph{$M$-free} if it is isomorphic to $M^{(I)}$ for some set $I$.

We shall use the following form of Walker's decomposition theorem.

\begin{lemma}[Walker \cite{Walker}]\label{lem:Walker}
Let $\eta$ be an infinite cardinal.  If a module $Y$ is a direct sum of $\eta$-generated modules, then every direct summand of $Y$ is a direct sum of $\eta$-generated modules.
\end{lemma}

The next elementary support-reduction lemma is used repeatedly.

\begin{lemma}\label{lem:support}
Let $A$ and $B$ be modules, let $\delta$ be an infinite cardinal, and suppose that $A$ is $\delta$-generated and isomorphic to a direct summand of $B^{(I)}$.  Then there is a subset $J\subseteq I$ with $\abs J\leq\delta$ such that $A$ is isomorphic to a direct summand of $B^{(J)}$.
\end{lemma}

\begin{proof}
Choose a split monomorphism $u:A\to B^{(I)}$ and a retraction $v:B^{(I)}\to A$ with $vu=1_A$.  Let $G$ be a generating set of $A$ with $\abs G\leq\delta$.  Every element $u(g)$ has finite support in $I$.  Hence
\[
  J=\bigcup_{g\in G}\operatorname{supp}(u(g))
\]
has cardinality at most $\delta$.  Since $u(A)\subseteq B^{(J)}$, the map $u$ factors through the canonical inclusion $B^{(J)}\hookrightarrow B^{(I)}$.  Restricting $v$ to $B^{(J)}$ still gives a left inverse of $u$.  Thus $A$ is a direct summand of $B^{(J)}$.
\end{proof}

We shall also use Eilenberg's classical swindle, in a form implicit in Bass's work on big projective modules \cite{Bass}; see also \cite{GobelTrlifaj}.

\begin{lemma}\label{lem:Eilenberg}
Let $A$ be a direct summand of a module $B$.  If $B\cong B^{(\aleph_0)}$, then
\[
  A\oplus B\cong B.
\]
\end{lemma}

\begin{proof}
Write $B\cong A\oplus C$.  Then
\[
 B\cong B^{(\aleph_0)}
 \cong A^{(\aleph_0)}\oplus C^{(\aleph_0)}.
\]
Consequently,
\[
 A\oplus B
 \cong A\oplus A^{(\aleph_0)}\oplus C^{(\aleph_0)}
 \cong A^{(\aleph_0)}\oplus C^{(\aleph_0)}
 \cong B.
\]
\end{proof}

\begin{lemma}\label{lem:stableSB}
Suppose that $A$ and $B$ are mutual direct summands and that
\[
 A\cong A^{(\aleph_0)},\qquad B\cong B^{(\aleph_0)}.
\]
Then $A\cong B$.
\end{lemma}

\begin{proof}
By \Cref{lem:Eilenberg}, $A\oplus B\cong B$ because $A$ is a direct summand of $B$, while $A\oplus B\cong A$ because $B$ is a direct summand of $A$.  Hence $A\cong B$.
\end{proof}

\section{The one-cardinal transfer theorem}

Fix from now on a nonzero module $M$ and an infinite cardinal $\eta$ such that $\gen(M)\leq\eta$.  Put
\[
  \theta=\max\{\abs R,\eta\}.
\]
Then $\abs M\leq\theta$ because $M$ is a quotient of $R^{(\eta)}$.

The first observation is the uniform decomposition bound that must be used in Walker's theorem.

\begin{proposition}\label{prop:small-decomp}
Every $X\in\Add(M)$ is a direct sum
\[
  X\cong\bigoplus_{i\in I}X_i
\]
where every $X_i$ is $\eta$-generated, belongs to $\Add(M)$, and satisfies $\abs{X_i}\leq\theta$.  The zero summands may be omitted.
\end{proposition}

\begin{proof}
By definition, $X$ is a direct summand of $M^{(J)}$ for some set $J$.  The module $M^{(J)}$ is a direct sum of $\eta$-generated modules, namely copies of $M$.  By Walker's theorem, $X$ is a direct sum of $\eta$-generated modules $X_i$.  Each $X_i$ is a direct summand of $X$, so $X_i\in\Add(M)$.  Finally, every $\eta$-generated module is a quotient of $R^{(\eta)}$, and therefore has cardinality at most $\theta$.
\end{proof}

The next theorem is the basic transfer statement.  It already shows that categoricity in a single sufficiently large cardinal implies categoricity in a tail, without any deconstructibility hypothesis.

\begin{theorem}\label{thm:upward}
Assume that $\Add(M)$ is $\lambda_0$-categorical for some cardinal $\lambda_0\geq\theta$.  Then, for every cardinal $\lambda\geq\lambda_0$, every $X\in\Add(M)$ with $\abs X=\lambda$ satisfies
\[
  X\cong M^{(\lambda)}.
\]
In particular,
\[
  [\lambda_0,\infty)\subseteq \Spec(\Add(M)).
\]
\end{theorem}

\begin{proof}
Let $X\in\Add(M)$ have cardinality $\lambda\geq\lambda_0$.  Since $M\neq0$ and $|M|\leq\theta\leq\lambda_0$, the module $M^{(\lambda_0)}$ has cardinality $\lambda_0$.  If $\lambda=\lambda_0$, then $X\cong M^{(\lambda_0)}$ by the assumed categoricity.  Hence assume $\lambda>\lambda_0\geq\theta$.  By \Cref{prop:small-decomp},
\[
  X\cong\bigoplus_{i\in I}X_i
\]
with $0\neq X_i$, $\abs{X_i}\leq\theta$.  We claim that $\abs I=\lambda$.  Since the summands are nonzero, $\abs I\leq\abs X=\lambda$.  Conversely, if $\abs I<\lambda$, then
\[
  \abs X\leq\max\{\abs I,\theta\}<\lambda,
\]
a contradiction.  Thus $\abs I=\lambda$.

Partition $I$ into pairwise disjoint subsets $I_\alpha$ of cardinality $\lambda_0$.  Put
\[
  Y_\alpha=\bigoplus_{i\in I_\alpha}X_i.
\]
Because every $X_i$ is nonzero, $|I_\alpha|=\lambda_0$, and $\abs{X_i}\leq\theta\leq\lambda_0$, one has $\abs{Y_\alpha}=\lambda_0$.  Moreover $Y_\alpha\in\Add(M)$.  Since $M\neq0$ and $\abs M\leq\theta\leq\lambda_0$, the module $M^{(\lambda_0)}$ also has cardinality $\lambda_0$.  By $\lambda_0$-categoricity,
\[
  Y_\alpha\cong M^{(\lambda_0)}
\]
for every $\alpha$.  Hence
\[
 X\cong\bigoplus_\alpha M^{(\lambda_0)}
 \cong M^{(\lambda)}.
\]
\end{proof}

\begin{remark}\label{rem:downward-gap}
The proof of \Cref{thm:upward} works only when $\lambda\geq\lambda_0$.  If $\theta<\lambda<\lambda_0$, the assumption that $\Add(M)$ is $\lambda_0$-categorical gives no direct information about two modules of cardinality $\lambda$.  Thus an assertion of the form ``$\abs P=\lambda$, $\abs{M^{(\lambda)}}=\lambda$, hence $P\cong M^{(\lambda)}$ by $\lambda_0$-categoricity'' is invalid.  A genuine downward transfer requires additional structure.  The next section identifies the relevant structure precisely.
\end{remark}

\section{Small-summand spectra and a stabilization criterion}

Let $\mathscr S_\eta(M)$ be a representative set of the isomorphism classes of all nonzero $\eta$-generated modules in $\Add(M)$.  Define
\[
  \sigma_M=\abs{\mathscr S_\eta(M)},
  \qquad
  \nu_M=\max\{\theta,\sigma_M\}.
\]
The first point is that $\mathscr S_\eta(M)$ is indeed a set of controlled cardinality.

\begin{lemma}\label{lem:number-types}
One has
\[
  \sigma_M\leq 2^\theta,
  \qquad\text{hence}\qquad
  \nu_M\leq 2^\theta.
\]
\end{lemma}

\begin{proof}
Every $\eta$-generated module is a quotient of $R^{(\eta)}$, and
\[
  \abs{R^{(\eta)}}=\theta.
\]
The set of submodules of $R^{(\eta)}$ has cardinality at most $2^\theta$.  Therefore the number of isomorphism classes of $\eta$-generated modules is at most $2^\theta$, and the same bound holds for the subclass $\mathscr S_\eta(M)$.
\end{proof}

Categoricity in one large cardinal forces all small summands to become indistinguishable after $\eta$-fold stabilization.

\begin{proposition}\label{prop:local-stab}
If $\Add(M)$ is $\lambda_0$-categorical for some $\lambda_0\geq\theta$, then
\[
  P^{(\eta)}\cong M^{(\eta)}
\]
for every $P\in\mathscr S_\eta(M)$.
\end{proposition}

\begin{proof}
Fix $P\in\mathscr S_\eta(M)$.  Since $P$ and $M$ are nonzero and have cardinality at most $\theta$, both $P^{(\lambda_0)}$ and $M^{(\lambda_0)}$ have cardinality $\lambda_0$.  They belong to $\Add(M)$, so
\[
  P^{(\lambda_0)}\cong M^{(\lambda_0)}.
\]

Because $P\in\Add(M)$ and $P$ is $\eta$-generated, \Cref{lem:support} shows that $P$ is a direct summand of $M^{(J)}$ for some $\abs J\leq\eta$.  Therefore $P^{(\eta)}$ is a direct summand of $M^{(\eta)}$.

Conversely, the displayed isomorphism shows that $M$ is a direct summand of $P^{(\lambda_0)}$.  Since $M$ is $\eta$-generated, another application of \Cref{lem:support} yields that $M$ is a direct summand of $P^{(J')}$ for some $\abs{J'}\leq\eta$.  Hence $M^{(\eta)}$ is a direct summand of $P^{(\eta)}$.

Now
\[
  \bigl(P^{(\eta)}\bigr)^{(\aleph_0)}\cong P^{(\eta)},
  \qquad
  \bigl(M^{(\eta)}\bigr)^{(\aleph_0)}\cong M^{(\eta)},
\]
because $\eta$ is infinite.  Thus \Cref{lem:stableSB} applies and gives $P^{(\eta)}\cong M^{(\eta)}$.
\end{proof}

We can now state the main theorem.

\begin{theorem}\label{thm:main}
Let $M\neq0$, let $\eta$ be an infinite cardinal with $\gen(M)\leq\eta$, and let $\theta$, $\mathscr S_\eta(M)$, and $\nu_M$ be as above.  The following conditions are equivalent.
\begin{enumerate}[label=\textup{(\arabic*)}]
\item $\Add(M)$ is $\lambda$-categorical for some cardinal $\lambda\geq\theta$.
\item For every $P\in\mathscr S_\eta(M)$,
\[
  P^{(\eta)}\cong M^{(\eta)}.
\]
\item For every nonzero $\eta$-generated module $P\in\Add(M)$,
\[
  \Add(P)=\Add(M).
\]
\item For every cardinal $\lambda>\nu_M$, every $X\in\Add(M)$ of cardinality $\lambda$ is $M$-free; more precisely,
\[
  X\cong M^{(\lambda)}.
\]
\item $\Add(M)$ is $\lambda$-categorical for every cardinal $\lambda>\nu_M$.
\end{enumerate}
\end{theorem}

\begin{proof}
$(1)\Rightarrow(2)$ is \Cref{prop:local-stab}.  We first prove $(2)\Longleftrightarrow(3)$.  Let $0\neq P\in\Add(M)$ be $\eta$-generated.  Since $P\in\Add(M)$, the support-reduction lemma gives
\[
   P\mid M^{(I)}
   \qquad\text{for some }|I|\leq\eta,
\]
where $A\mid B$ means that $A$ is isomorphic to a direct summand of $B$.  Hence $P^{(\eta)}\mid M^{(\eta)}$.

If (2) holds, choose the representative $P_0\in\mathscr S_\eta(M)$ with $P\cong P_0$.  Then $P^{(\eta)}\cong P_0^{(\eta)}\cong M^{(\eta)}$, so in particular $M\in\Add(P)$.  Since already $P\in\Add(M)$, closure under direct sums and direct summands yields
\[
   \Add(P)=\Add(M).
\]
Thus (2) implies (3).

Conversely, assume (3).  Then $M\in\Add(P)$.  As $M$ is $\eta$-generated, support reduction gives $M\mid P^{(J)}$ for some $|J|\leq\eta$, and therefore $M^{(\eta)}\mid P^{(\eta)}$.  We already know $P^{(\eta)}\mid M^{(\eta)}$.  Both modules are isomorphic to their countable direct sums because $\eta$ is infinite.  The stable Schr\"oder--Bernstein lemma therefore gives
\[
  P^{(\eta)}\cong M^{(\eta)}.
\]
Hence (3) implies (2).

The implication $(4)\Rightarrow(5)$ is immediate, and $(5)\Rightarrow(1)$ follows by choosing any $\lambda>\nu_M\geq\theta$.  It remains to prove $(2)\Rightarrow(4)$.

Assume (2), let $\lambda>\nu_M$, and let $X\in\Add(M)$ with $|X|=\lambda$.  By \Cref{prop:small-decomp}, and after grouping isomorphic summands, we can write
\[
   X\cong\bigoplus_{P\in\mathscr S_\eta(M)}P^{(\kappa_P)}
\]
for suitable cardinals $\kappa_P$.  Put
\[
   \mathscr T=\{P\in\mathscr S_\eta(M)\mid \kappa_P<\eta\},
\]
and define
\[
   X_{<}=\bigoplus_{P\in\mathscr T}P^{(\kappa_P)},
   \qquad
   X_{\geq}=\bigoplus_{P\notin\mathscr T}P^{(\kappa_P)}.
\]
There are at most $\sigma_M\eta\leq\nu_M$ summands in $X_{<}$, and every one has cardinality at most $\theta\leq\nu_M$.  Hence
\[
   |X_{<}|\leq\nu_M<\lambda.
\]

For $P\notin\mathscr T$ we have $\kappa_P\geq\eta$, and therefore
\[
\begin{aligned}
 P^{(\kappa_P)}
   &\cong \bigl(P^{(\eta)}\bigr)^{(\kappa_P)}
    \cong \bigl(M^{(\eta)}\bigr)^{(\kappa_P)}
    \cong M^{(\kappa_P)}.
\end{aligned}
\]
Consequently
\[
  X_{\geq}\cong M^{(\mu)},
  \qquad
  \mu=\sum_{P\notin\mathscr T}\kappa_P.
\]
Since $X=X_<\oplus X_{\geq}$ has cardinality $\lambda$ whereas $|X_<|<\lambda$, we have $|X_{\geq}|=\lambda$.  As $M\neq0$ and $|M|\leq\theta<\lambda$, it follows that $\mu=\lambda$.  Thus
\[
   X_{\geq}\cong M^{(\lambda)}.
\]

It remains only to absorb $X_<$.  Since $X_<\in\Add(M)$ and
\[
   \gen(X_<)\leq |X_<|\leq\nu_M,
\]
the support-reduction lemma, applied with $\delta=\nu_M$, shows that $X_<$ is a direct summand of $M^{(J)}$ for some $|J|\leq\nu_M<\lambda$.  Hence $X_<$ is a direct summand of $M^{(\lambda)}$.  Since
\[
   M^{(\lambda)}\cong\bigl(M^{(\lambda)}\bigr)^{(\aleph_0)},
\]
Eilenberg absorption gives
\[
   X_<\oplus M^{(\lambda)}\cong M^{(\lambda)}.
\]
Therefore
\[
  X\cong X_<\oplus X_{\geq}
    \cong X_<\oplus M^{(\lambda)}
    \cong M^{(\lambda)}.
\]
This proves (4) and completes the proof.
\end{proof}

The theorem has two useful immediate consequences.

\begin{corollary}\label{cor:2theta}
If $\Add(M)$ is categorical in some cardinal $\lambda\geq\theta$, then it is categorical in every cardinal
\[
  \kappa>2^\theta,
\]
and every $X\in\Add(M)$ of cardinality $\kappa>2^\theta$ is isomorphic to $M^{(\kappa)}$.
\end{corollary}

\begin{proof}
By \Cref{lem:number-types}, $\nu_M\leq2^\theta$.  Apply \Cref{thm:main}.
\end{proof}

\begin{corollary}\label{cor:smalltypes}
Assume in addition that
\[
  \abs{\mathscr S_\eta(M)}\leq\theta.
\]
Then the following are equivalent.
\begin{enumerate}[label=\textup{(\arabic*)}]
\item $\Add(M)$ is categorical in some cardinal $\lambda\geq\theta$.
\item $\Add(M)$ is categorical in every cardinal $\lambda>\theta$.
\item Every $X\in\Add(M)$ of cardinality $\lambda>\theta$ satisfies $X\cong M^{(\lambda)}$.
\item Every $X\in\Add(M)$ of cardinality $\theta^+$ satisfies $X\cong M^{(\theta^+)}$.
\end{enumerate}
\end{corollary}

\begin{proof}
The hypothesis gives $\nu_M=\theta$, so $(1)\Rightarrow(3)\Rightarrow(2)$ follows from \Cref{thm:main}.  Clearly $(3)\Rightarrow(4)$.  Finally, (4) implies that $\Add(M)$ is $\theta^+$-categorical because $M^{(\theta^+)}$ has cardinality $\theta^+$, hence (1) follows.
\end{proof}

\begin{remark}\label{rem:missing-parameter}
The role of $\sigma_M$ is essential in the proof.  Walker's theorem bounds the size of each small summand, but it does not bound the number of their isomorphism types by $\theta$; a priori that number can be as large as $2^\theta$.  Thus a transfer theorem with threshold $\theta$ requires either an argument that controls these types or an explicit hypothesis such as the one in \Cref{cor:smalltypes}.  The unconditional statement furnished by \Cref{thm:main} is the tail above $\nu_M$, and \Cref{cor:2theta} gives a bound depending only on $R$ and the number of generators of $M$.
\end{remark}

Theorems on deconstructible classes show that sufficiently high categoricity often propagates to a tail; see \cite{Trlifaj,SarochTrlifaj}.  The mechanism here is different.  The class $\Add(M)$ need not be extension-closed, and the proof uses Walker decomposition, stabilization of small summands, and Eilenberg absorption rather than filtrations.

\begin{corollary}\label{cor:spectrum}
	With the notation of \Cref{thm:main}, suppose
	\[
	\Spec(\Add(M))\cap[\theta,\infty)\neq\varnothing.
	\]
	Then:
	\begin{enumerate}[label=\textup{(\roman*)}]
		\item for every $\lambda_0\in\Spec(\Add(M))$ with $\lambda_0\geq\theta$,
		\[
		[\lambda_0,\infty)\subseteq\Spec(\Add(M));
		\]
		\item
		\[
		(\nu_M,\infty)\subseteq\Spec(\Add(M));
		\]
		\item in particular,
		\[
		(2^\theta,\infty)\subseteq\Spec(\Add(M)).
		\]
	\end{enumerate}
\end{corollary}

\begin{proof}
	Part (i) is the upward-transfer theorem; part (ii) is \Cref{thm:main}; and part (iii) follows from \Cref{lem:number-types}.
\end{proof}

\begin{proposition}\label{prop:bounded-remainder}
	Assume the equivalent conditions of \Cref{thm:main}.  For every $X\in\Add(M)$ of cardinality $\lambda>\nu_M$, there is a decomposition
	\[
	X\cong X_0\oplus M^{(\lambda)}
	\]
	with $|X_0|\leq\nu_M$.  Moreover $X_0$ is a direct summand of $M^{(\lambda)}$, and hence $X\cong M^{(\lambda)}$.
\end{proposition}

\begin{proof}
	Take $X_0=X_<$ in the proof of \Cref{thm:main}.  The cardinality bound, the decomposition, and the fact that $X_0$ is a direct summand of $M^{(\lambda)}$ were all established there.
\end{proof}

The proposition shows that the obstruction to immediate $M$-freeness is concentrated in a bounded part formed from small isomorphism types occurring with multiplicity $<\eta$.  Once the high-multiplicity part reaches the ambient cardinality, the bounded remainder is absorbed.

\begin{corollary}\label{cor:smalltype-restated}
	If $|\mathscr S_\eta(M)|\leq\theta$ and $\Add(M)$ is categorical in some cardinal $\lambda\geq\theta$, then every $X\in\Add(M)$ of cardinality $\kappa>\theta$ satisfies
	\[
	X\cong M^{(\kappa)}.
	\]
	Thus $\Add(M)$ is categorical in every cardinal $\kappa>\theta$.
\end{corollary}

\begin{proof}
	Here $\nu_M=\theta$, so the assertion is immediate from \Cref{thm:main}.
\end{proof}

\section{Projective modules, trace ideals, and \texorpdfstring{$p$}{p}-connected rings}\label{sec:projective}

Assume throughout this section that $R\neq0$.  Put
\[
   \theta_R=\max\{|R|,\aleph_0\}.
\]
Let $\mathscr P_c(R)$ be a representative set of the nonzero countably generated projective right $R$-modules, and put
\[
   \nu_R=\max\{\theta_R,|\mathscr P_c(R)|\}.
\]
By Kaplansky's theorem every projective module is a direct sum of countably generated projective modules \cite{Kaplansky}.  In particular, the general theorem with $M=R$ and $\eta=\aleph_0$ recovers the stabilization criterion of Trlifaj \cite[Proposition~1.1]{Trlifaj}.

For a right $R$-module $P$, its \emph{trace ideal} is
\[
  \operatorname{Tr}_R(P)
   =\sum_{f\in\operatorname{Hom}_R(P,R)}f(P),
\]
a two-sided ideal of $R$.  A projective module $P$ is a generator of $\operatorname{Mod}\text{-}R$ if and only if $\operatorname{Tr}_R(P)=R$; equivalently, $R$ is a direct summand of a finite direct sum of copies of $P$ \cite{AndersonFuller}.  Following Bass \cite{Bass}, we call $R$ \emph{right $p$-connected} if every nonzero projective right $R$-module is a generator.  Bass proved, in particular, that over a $p$-connected ring the countable direct sum of copies of any projective module is free \cite[Proposition~2.4]{Bass}.  We include a short proof of the precise stabilization form needed here.

\begin{lemma}\label{lem:generator-stab}
Let $0\neq P$ be a countably generated projective right $R$-module.  If $P$ is a generator, then
\[
   P^{(\aleph_0)}\cong R^{(\aleph_0)}.
\]
\end{lemma}

\begin{proof}
Because $P$ is countably generated and projective, $P$ is a direct summand of $R^{(\aleph_0)}$, whence
\[
   P^{(\aleph_0)}\mid R^{(\aleph_0)}.
\]
Since $P$ is a generator, there is an integer $n\geq1$ such that $R\mid P^n$.  Taking countable direct sums gives
\[
   R^{(\aleph_0)}\mid P^{(\aleph_0)}.
\]
Both modules are isomorphic to their countable direct sums, so \Cref{lem:stableSB} yields the desired isomorphism.
\end{proof}

\begin{theorem}\label{thm:pconnected}
The following conditions are equivalent.
\begin{enumerate}[label=\textup{(\arabic*)}]
\item $\Proj$ is categorical in some cardinal $\lambda\geq\theta_R$.
\item For every nonzero countably generated projective module $P$,
\[
   P^{(\aleph_0)}\cong R^{(\aleph_0)}.
\]
\item $R$ is right $p$-connected.
\item Every projective module of cardinality $\lambda>\nu_R$ is free; more precisely, it is isomorphic to $R^{(\lambda)}$.
\item $\Proj$ is categorical in every cardinal $\lambda>\nu_R$.
\end{enumerate}
Moreover $\nu_R\leq2^{\theta_R}$.
\end{theorem}

\begin{proof}
The equivalence of (1), (2), (4), and (5), as well as the bound $\nu_R\leq2^{\theta_R}$, follows from \Cref{thm:main,lem:number-types} with $M=R$ and $\eta=\aleph_0$.

Assume (3).  Every nonzero countably generated projective module is a generator, so (2) follows from \Cref{lem:generator-stab}.

Conversely, assume (2), and let $0\neq Q$ be an arbitrary projective module.  By Kaplansky's theorem,
\[
   Q\cong\bigoplus_{i\in I}Q_i
\]
with each $Q_i$ countably generated and projective.  Choose a nonzero summand $P=Q_{i_0}$.  By (2), $R$ is a direct summand of $P^{(\aleph_0)}$.  Choose a split monomorphism $u:R\to P^{(\aleph_0)}$.  The element $u(1)$ has finite support, so $u$ factors through $P^n$ for some finite $n$; restricting a retraction of $u$ shows
\[
   R\mid P^n.
\]  Thus $P$ is a generator.  Since $P$ is a direct summand of $Q$, the trace ideal of $Q$ contains that of $P$, and therefore $Q$ is a generator.  Hence $R$ is right $p$-connected.
\end{proof}

\begin{corollary}\label{cor:containsR}
Let $M\neq0$ be $\eta$-generated for an infinite cardinal $\eta$, assume $R\in\Add(M)$, and put $\theta=\max\{|R|,\eta\}$.  The following are equivalent.
\begin{enumerate}[label=\textup{(\arabic*)}]
\item $\Add(M)$ is categorical in some cardinal $\lambda\geq\theta$.
\item $M$ is projective and $R$ is right $p$-connected.
\item $\Add(M)=\Proj$ and $\Proj$ is categorical in every cardinal $\lambda>\nu_R$.
\end{enumerate}
When these conditions hold, every module in $\Add(M)$ of cardinality $\lambda>\nu_R$ is free.
\end{corollary}

\begin{proof}
Assume (1), and let $\lambda\geq\theta$ be a categoricity cardinal.  Since $R\in\Add(M)$, both $R^{(\lambda)}$ and $M^{(\lambda)}$ belong to $\Add(M)$; both have cardinality $\lambda$.  Hence
\[
   M^{(\lambda)}\cong R^{(\lambda)}.
\]
The module $M$ is a direct summand of the left-hand side, and therefore is projective.  Thus $\Add(M)\subseteq\Proj$.  On the other hand, $R\in\Add(M)$ and closure under direct sums and summands give
\[
   \Proj=\Add(R)\subseteq\Add(M).
\]
So $\Add(M)=\Proj$.  Since this class is categorical in a cardinal $\lambda\geq\theta\geq\theta_R$, \Cref{thm:pconnected} implies that $R$ is right $p$-connected.  This proves (2).

If (2) holds, then $M$ projective gives $\Add(M)\subseteq\Proj$, while $R\in\Add(M)$ gives the reverse inclusion.  Thus $\Add(M)=\Proj$, and \Cref{thm:pconnected} gives (3).  Finally (3) implies (1) by choosing a cardinal larger than both $\theta$ and $\nu_R$.
\end{proof}

\begin{remark}\label{rem:Mfree-free}
The hypothesis $R\mid M$ alone does not identify $M$-free modules with free modules.  For example, over $R=\mathbb Z$ let
\[
   M=\mathbb Z\oplus\mathbb Z/2\mathbb Z.
\]
Then $R$ is a direct summand of $M$, but $M^{(I)}$ has nonzero $2$-torsion for every nonempty set $I$, and hence is not a free abelian group.  What categoricity supplies in \Cref{cor:containsR} is precisely the missing projectivity of $M$, after which the whole class reduces to $\Proj$.
\end{remark}

\section{Semisimple modules}

The semisimple case is simpler and admits a direct proof with the optimal natural threshold.  This avoids any need to form a single generator from all isomorphism classes of simple modules; the number of such classes can be larger than $\abs R$.

\begin{theorem}\label{thm:semisimple}
Let $R$ be nonzero.  Let $\SSmod$ denote the class of all semisimple right $R$-modules and let
\[
  \theta_R=\max\{\abs R,\aleph_0\}.
\]
The following are equivalent.
\begin{enumerate}[label=\textup{(\arabic*)}]
\item $\SSmod$ is categorical in some cardinal $\lambda>\theta_R$.
\item $\SSmod$ is categorical in every cardinal $\lambda>\theta_R$.
\item $R$ admits a unique simple right module up to isomorphism.
\end{enumerate}
If these conditions hold and $S$ is the unique simple module, then every semisimple module $X$ of cardinality $\lambda>\theta_R$ satisfies
\[
  X\cong S^{(\lambda)}.
\]
\end{theorem}

\begin{proof}
Assume first that $R$ has a unique simple right module $S$ up to isomorphism.  Every semisimple module is isomorphic to $S^{(\kappa)}$ for some cardinal $\kappa$.  Since $S$ is cyclic, $\abs S\leq\abs R\leq\theta_R$.  If $\abs{S^{(\kappa)}}=\lambda>\theta_R$, then necessarily $\kappa=\lambda$.  Thus every semisimple module of cardinality $\lambda$ is isomorphic to $S^{(\lambda)}$, proving (2).

Clearly (2) implies (1).  Conversely, suppose that $R$ has two nonisomorphic simple right modules $S$ and $T$.  For every $\lambda>\theta_R$, both $S^{(\lambda)}$ and $T^{(\lambda)}$ have cardinality $\lambda$.  They cannot be isomorphic.  Indeed, since $S$ is cyclic, every map $S\to T^{(\lambda)}$ has image contained in a finite subsum; but $\operatorname{Hom}_R(S,T)=0$, so every such map is zero.  Hence $S^{(\lambda)}\not\cong T^{(\lambda)}$.  Therefore $\SSmod$ is not categorical in any cardinal $\lambda>\theta_R$, proving (1)$\Rightarrow$(3).
\end{proof}

\begin{example}\label{ex:local-semisimple}
Let $R=\mathbb Z/4\mathbb Z$.  Up to isomorphism, $R$ has a unique simple right module, namely $S=\mathbb Z/2\mathbb Z$.  Hence \Cref{thm:semisimple} shows that the class of semisimple $R$-modules is categorical in every cardinal $\lambda>\aleph_0$.  Nevertheless $R$ is not a semisimple ring.  Thus categoricity of the \emph{class of semisimple modules} detects uniqueness of the simple isomorphism type, not semisimplicity of the ambient ring.
\end{example}

\begin{remark}\label{rem:semisimple-cardinality}
One may write $\SSmod=\Add(S_0)$ where $S_0$ is a direct sum of representatives of all simple modules.  However, the number of simple isomorphism types is bounded a priori only by $2^{\abs R}$, not by $\abs R$.  Thus deriving \Cref{thm:semisimple} from a theorem whose generating bound is $\abs R+\aleph_0$ requires an additional cardinality argument.  The direct proof above is both shorter and sharper.
\end{remark}

\section{Pure-projective modules and simple von Neumann regular rings}\label{sec:pureproj}

Let $\PProj$ denote the class of all pure-projective right $R$-modules.  We use the classical characterization that a module is pure-projective if and only if it is a direct summand of a direct sum of finitely presented modules \cite{GobelTrlifaj,Warfield}.  Walker's theorem immediately gives the following decomposition.

\begin{lemma}\label{lem:pp-countable}
Every pure-projective module is a direct sum of countably generated pure-projective modules.
\end{lemma}

\begin{proof}
A pure-projective module is a direct summand of a direct sum of finitely presented modules.  Every finitely presented module is finitely generated, hence countably generated.  Walker's theorem therefore decomposes the direct summand into countably generated modules, and each of those summands is again pure-projective.
\end{proof}

We first identify the interaction between Bass's connectedness condition and von Neumann regularity.

\begin{proposition}\label{prop:simple-vnr-pconn}
Let $R\neq0$ be a von Neumann regular ring.  Then the following are equivalent.
\begin{enumerate}[label=\textup{(\arabic*)}]
\item $R$ is right $p$-connected.
\item $R$ is simple.
\end{enumerate}
\end{proposition}

\begin{proof}
Assume first that $R$ is simple.  Let $0\neq P$ be projective.  We claim that $\operatorname{Tr}_R(P)\neq0$.  Choose $0\neq p\in P$.  Since $P$ is projective, the dual-basis lemma gives elements $p_1,\dots,p_n\in P$ and maps $f_1,\dots,f_n\in\operatorname{Hom}_R(P,R)$ such that
\[
   p=\sum_{i=1}^n p_i f_i(p).
\]
Thus $f_i(p)\neq0$ for at least one $i$, and so the trace ideal is nonzero.  Being a nonzero two-sided ideal in a simple ring, it equals $R$.  Hence every nonzero projective module is a generator and $R$ is right $p$-connected.

Conversely, assume that $R$ is right $p$-connected, and let $0\neq I\triangleleft R$ be a two-sided ideal.  Choose $0\neq a\in I$.  Since $R$ is von Neumann regular, there is $x\in R$ with $a=axa$.  Put $e=ax$.  Then $e^2=e$, $e\neq0$, and $e\in I$.  The right ideal $eR$ is a nonzero finitely generated projective module.  By $p$-connectedness,
\[
   \operatorname{Tr}_R(eR)=R.
\]
For an idempotent $e$, the trace ideal of $eR$ is $ReR$.  Since $e\in I$ and $I$ is two-sided,
\[
   R=ReR\subseteq I.
\]
Thus $I=R$, proving that $R$ is simple.
\end{proof}

The main result of this section is now completely ring-theoretic.

\begin{theorem}\label{thm:pureproj}
Assume $R\neq0$, and put
\[
   \theta_R=\max\{|R|,\aleph_0\},
   \qquad
   \nu_R=\max\{\theta_R,|\mathscr P_c(R)|\}.
\]
The following conditions are equivalent.
\begin{enumerate}[label=\textup{(\arabic*)}]
\item $\PProj$ is categorical in some cardinal $\lambda\geq\theta_R$.
\item For every nonzero countably generated pure-projective module $P$,
\[
   P^{(\aleph_0)}\cong R^{(\aleph_0)}.
\]
\item $R$ is von Neumann regular and right $p$-connected.
\item $R$ is a simple von Neumann regular ring.
\item $\PProj=\Proj$, and this class is categorical in every cardinal $\lambda>\nu_R$.
\item Every pure-projective module of cardinality $\lambda>\nu_R$ is free; more precisely, it is isomorphic to $R^{(\lambda)}$.
\end{enumerate}
In particular, if one of these conditions holds, then $\PProj$ is categorical in every cardinal $\lambda>2^{\theta_R}$.
\end{theorem}

\begin{proof}
$(1)\Rightarrow(2)$.  Let $0\neq P$ be countably generated and pure-projective, and let $\lambda\geq\theta_R$ be a categoricity cardinal.  Both $P^{(\lambda)}$ and $R^{(\lambda)}$ are pure-projective of cardinality $\lambda$, hence
\[
   P^{(\lambda)}\cong R^{(\lambda)}.
\]
Since $P$ is countably generated, support reduction shows that $P$ is a direct summand of $R^{(\aleph_0)}$; hence $P^{(\aleph_0)}\mid R^{(\aleph_0)}$.  Conversely, the displayed isomorphism makes $R$ a direct summand of $P^{(\lambda)}$, and because $R$ is cyclic, support reduction gives $R\mid P^{(\aleph_0)}$.  Thus $R^{(\aleph_0)}\mid P^{(\aleph_0)}$.  Stable Schr\"oder--Bernstein yields
\[
   P^{(\aleph_0)}\cong R^{(\aleph_0)}.
\]

$(2)\Rightarrow(3)$.  Let $F$ be a finitely presented right $R$-module.  Then $F$ is pure-projective.  If $F\neq0$, condition (2) makes $F$ a direct summand of the free module $F^{(\aleph_0)}\cong R^{(\aleph_0)}$, so $F$ is projective.  The zero module is projective as well.  Hence every finitely presented right module is projective, which is a standard characterization of von Neumann regular rings \cite{Goodearl,Lam}.  Thus $R$ is von Neumann regular.

Every countably generated projective module is pure-projective, so (2) also gives the stabilization condition in \Cref{thm:pconnected}; consequently $R$ is right $p$-connected.  This proves (3).

$(3)\Longleftrightarrow(4)$ is \Cref{prop:simple-vnr-pconn}.

Assume (3).  Over a von Neumann regular ring every finitely presented module is projective.  Hence every pure-projective module, being a direct summand of a direct sum of finitely presented modules, is projective; the reverse inclusion holds over every ring.  Thus
\[
   \PProj=\Proj.
\]
Right $p$-connectedness and \Cref{thm:pconnected} then imply (5).

$(5)\Rightarrow(6)$ follows because $R^{(\lambda)}$ is a pure-projective module of cardinality $\lambda$ for every $\lambda>\nu_R$.  Finally, $(6)\Rightarrow(1)$ by choosing any $\lambda>\nu_R$.  The final bound follows from $\nu_R\leq2^{\theta_R}$ in \Cref{thm:pconnected}.
\end{proof}

\begin{corollary}\label{cor:comm-pureproj}
Let $R\neq0$ be commutative.  The class of pure-projective $R$-modules is categorical in some cardinal $\lambda\geq\max\{|R|,\aleph_0\}$ if and only if $R$ is a field.
\end{corollary}

\begin{proof}
A commutative simple ring is a field.  Apply \Cref{thm:pureproj}.
\end{proof}

\begin{example}\label{ex:matrix}
Let $R=M_n(D)$ for a division ring $D$ and $n\geq2$.  Then $R$ is simple von Neumann regular, so $\PProj$ is eventually categorical.  This example shows why stabilization, rather than finite-rank freeness, is the natural phenomenon.  If $e$ is a primitive idempotent, then $eR$ is a nonzero finitely generated projective module which is not free as an $R$-module when $n>1$, but
\[
   (eR)^{(\aleph_0)}\cong R^{(\aleph_0)}.
\]
Indeed, $eR$ is a generator because $R$ is simple, and \Cref{lem:generator-stab} applies.
\end{example}

\begin{example}\label{ex:ultramatricial}
Let $K$ be a field and consider the unital direct system
\[
 M_1(K)\longrightarrow M_2(K)\longrightarrow M_4(K)\longrightarrow\cdots,
 \qquad
 a\longmapsto\begin{pmatrix}a&0\\0&a\end{pmatrix}.
\]
Its direct limit $R$ is a simple von Neumann regular ring; simplicity of this standard direct-limit example is also used in \cite[Example~2.4]{Trlifaj}.  Von Neumann regularity follows because every element lies in a matrix stage and already has an inner inverse there.  The ring is not artinian; for example, it contains infinitely many nonzero pairwise orthogonal idempotents.  Hence \Cref{thm:pureproj} produces eventual categoricity of pure-projective modules far beyond the semisimple-artinian setting.
\end{example}


\begin{thebibliography}{99}

\bibitem{AndersonFuller}
F. W. Anderson and K. R. Fuller,
\emph{Rings and Categories of Modules}, 2nd ed., Graduate Texts in Mathematics 13, Springer, New York, 1992.

\bibitem{Bass}
H. Bass,
Big projective modules are free,
\emph{Illinois J. Math.} \textbf{7} (1963), 24--31.

\bibitem{GobelTrlifaj}
R. G\"obel and J. Trlifaj,
\emph{Approximations and Endomorphism Algebras of Modules}, 2nd ed., De Gruyter Expositions in Mathematics 41, Vols. 1--2, De Gruyter, Berlin--Boston, 2012.

\bibitem{Goodearl}
K. R. Goodearl,
\emph{Von Neumann Regular Rings}, 2nd ed., Krieger Publishing, Malabar, FL, 1991.

\bibitem{Kaplansky}
I. Kaplansky,
Projective modules,
\emph{Ann. of Math. (2)} \textbf{68} (1958), 372--377.

\bibitem{Lam}
T. Y. Lam,
\emph{Lectures on Modules and Rings}, Graduate Texts in Mathematics 189, Springer, New York, 1999.

\bibitem{Los}
J. \L o\'s,
On the categoricity in power of elementary deductive systems and some related problems,
\emph{Colloq. Math.} \textbf{3} (1954), 58--62.

\bibitem{Mazari}
M. Mazari-Armida,
Characterizing categoricity in several classes of modules,
\emph{J. Algebra} \textbf{617} (2023), 382--401.

\bibitem{Morley}
M. Morley,
Categoricity in power,
\emph{Trans. Amer. Math. Soc.} \textbf{114} (1965), 514--538.

\bibitem{SarochTrlifaj}
J. \v{S}aroch and J. Trlifaj,
Deconstructible abstract elementary classes of modules and categoricity,
\emph{Bull. Lond. Math. Soc.} \textbf{56} (2024), 3854--3866.

\bibitem{Shelah}
S. Shelah,
On what I do not understand (and have something to say), model theory,
\emph{Math. Japon.} \textbf{51} (2000), 329--377.

\bibitem{Trlifaj}
J. Trlifaj,
Categoricity for transfinite extensions of modules,
\emph{Proc. Amer. Math. Soc. Ser. B} \textbf{10} (2023), 369--381.

\bibitem{Walker}
C. P. Walker,
Relative homological algebra and Abelian groups,
\emph{Illinois J. Math.} \textbf{10} (1966), 186--209.

\bibitem{Warfield}
R. B. Warfield, Jr.,
Purity and algebraic compactness for modules,
\emph{Pacific J. Math.} \textbf{28} (1969), 699--719.

\end{thebibliography}
\end{document}